\documentclass[a4paper,12pt,legno]{article}
\usepackage[polish]{babel}
\usepackage[cp1250]{inputenc}
\usepackage[T1]{fontenc}
\usepackage{amsmath}
\usepackage{amsfonts}
\usepackage{amsthm}
\usepackage{mathrsfs}

\newtheorem{stw}[]{Proposition}
\newtheorem{lem}[]{Lemma}

\newtheorem{wn}{Corollary}

\newtheorem{rem}[]{Remark}

\newcommand{\Rset}{\mathbb{R}} 
\newcommand{\J}{\mathcal{J}}   
\newcommand{\La}{\mathcal{L}}  

\newcommand{\nn}{\nonumber}
\newcommand{\al}{\alpha}
\newcommand{\be}{\beta}

\newcommand{\bl}{\bigl(}
\newcommand{\br}{\bigr)}
\newcommand{\Bl}{\Bigl(}
\newcommand{\Br}{\Bigr)}
\newcommand{\wh}{^{\widehat{}}}

\begin{document}

\title{A note on a composition of two random integral mappings $\J^\be$ and some examples}

\author{Agnieszka Czyżewska-Jankowska and Zbigniew J. Jurek\footnote{Research funded by a grant MEN Nr 1P03A04629,
2005-2008.}}

\date{November 15, 2008}
\maketitle
\begin{quote} \textbf{Abstract.} A method of \emph{random integral
representation}, that is, a method of representing a given
probability measure as the probability distribution of some random
integral, was quite successful in the past few decades. In this
note we show that a composition of two  random integral mappings
$\J^\be$ is again a random integral mapping. We illustrate our
results on some examples.



\emph{Abbreviated title:}  A Composition of random mappings $\J^\be$

\end{quote}

We say that a probability distribution (measure) $\mu$ admits
\emph{a random integral representation} if we have
 \emph{
\begin{multline}
\qquad \qquad \qquad \qquad \qquad \mu=\mathcal{L}\big( \int_I h(t)dY(r(t))\big),   \\
\mbox{where} \ I=(a,b] \subset \Rset^+, \
h:\Rset^+\to \Rset, \ \ Y(\cdot) \ \mbox{is a L\'evy process and}, \\
r:\Rset^+ \to \Rset^+  \mbox{is a monotone function (a time change
in $Y$)}
\end{multline}}
In fact, in the past it was proved that many classes of limit laws
can be described as probability distributions of random integrals of
the form (1). Moreover, it was  conjectured that \emph{all classes
of limit laws derived for sequences of independent random variables
} should admit a random integral representation; cf. Jurek (1985;
1988) and see the Conjecture on

\noindent $www.math.uni.wroc.pl/$$\sim$zjjurek  \ . \ The random
integral approach was also successfully  used by others; see for
instance Aoyama-Maejima (2007). Last but not least, one should
emphasis that from the integral representations (1) very easily
follow formulae for the characteristic functions and the L\'evy-
Khintchine representation.

\medskip
In this note we examine how some random integral mappings of the
form (1) behave under a composition of such mappings.

\medskip
\textbf{1. Introduction and main results.} Throughout the paper
$\mathcal{L}(X)$ will denote the probability distribution of a
$\Rset^d$-valued\footnote{Our proofs are a such that they hold true
for a real separable Banach space valued random elements as well.
Interested Readers in probability on Banach spaces, e.g., those who
prefer to see a stochastic process with continuous paths on $[0,1]$
interval as an $C[0,1]$-valued random element, we refer to the
monograph by Araujo and Gine (1980).} random vector $X$. Similarly,
by $Y_{\nu}(t), t\ge0,$ we will denote an $\Rset^d$- valued L\'evy
process such that $\mathcal{L}_{\nu}(Y(1))=\nu$. By a L\'evy process
we mean process starting from zero, with stationary and independent
increments and with paths that are right continuous and have finite
left hand limits. Of course, we always have that $\nu\in ID$, where
$ID$ stands for the set of all \emph{\underline{i}nfinitely
\underline{d}ivisible} measures on $\Rset^d$, so in particular
$\nu^{\ast c}$, $c>0$, is well-defined infinitely divisible
probability measure.

For $\beta>0$ and a L\'evy process $Y_{\nu}(t), t\ge0$, we define an
integral mapping  \qquad $ \mathcal{J}^{\beta}:\, ID\to ID$  and a
class $\mathcal{U}_{\beta}$ as follows
\begin{equation}
 \mathcal{J}^{\beta}(\nu):\,=\La\bl\int_0^1
t^{1/\beta}\;dY_{\nu}(t)\br=\La\bl\int_0^1 t\;dY_{\nu}(t^{\beta})
\br, \ \ \mbox{and}\ \
\mathcal{U}_{\beta}:\,=\mathcal{J}^{\beta}(ID).
\end{equation}
For another equivalent characterizations of classes
$\mathcal{U}_{\beta}$, even in a greater generality than we are
interested in the present note, we refer to Jurek (1988). To the
distributions from the class $\mathcal{U}_{\beta}$ we refer to as
\emph{generalized s-selfdecomposable distributions}.

Recall here that the importance and the interest in the increasing
family $\mathcal{U}_\be$ of convolution semigroups stems from the
fact that
\[
\overline{(\cup_n\mathcal{U}_{\be_n})}=ID, \ \ \mbox{for any
increasing sequeunce} \ \ \be_n \to \infty\, ,
\]
where the bar means a closure in the weak topology; cf. Jurek
(1996), Corollary 1, and the references therein. More explicitly,
let us note that
\[
\J^{\be}(\nu)\Rightarrow \mathcal{L}(\int_0^1\,dY_{\nu}(t)) =\nu, \
\ \mbox{as} \ \be \to \infty\,.
\]

\medskip
Since we assume that the paths of $Y$ are almost surely cadlag
(i.e., right continuous with finite left hand limits) and the random
integral in (2) we define by a formal integration by parts formula,
therefore the random integral in question exists; for details cf.
Jurek-Vervaat (1983), Lemma 1.1 or Jurek-Mason (1993), Section 3.6,
p. 116.

\medskip
Here are the results:
\begin{stw}
For $\nu \in ID$  and for positive $\al\neq \be$ we have
\begin{equation}
\J^\al(\J^\be(\nu)) =
\La(\int^{1}_{0}u\,dY_{\nu}(r_{(\al,\be)}(u)\br\bigr)=\La(\int^{1}_{0}\,
r_{(\al,\be)}^{-1}(u)\,dY_{\nu}(u)\br\bigr)
\end{equation}
where $r_{(\al,\be)} :[0,1]\to [0,1]$ is a continuous strictly
increasing time change given by the formula
$r_{(\al,\be)}(u):\,=\tfrac{\beta}{\beta-\alpha}u^{\al} -
\tfrac{\alpha}{\beta-\alpha}u^\be$, which is symmetric  in  $\al$
and $\be$, that is, $r_{(\al , \be)}(u)= r_{(\be,\al)}(u)$.
\end{stw}

\begin{wn}
 For  $0<\al<\be$  and $\nu \in ID$ we have the identity:
 \begin{equation}
\J^\alpha\big(\J^\beta(\nu^{\ast(\be-\al)})\bigr)*\J^\beta
(\nu^{\ast\alpha})= \J^\alpha\bigl(\nu^{\ast \beta}\bigr)
\end{equation}
Equivalently,  in terms of characteristic functions we have
\begin{equation}
(\be-\al)\log (\J^\alpha(\J^\beta(\nu)))\wh (y)=
\be\log(\J^\al(\nu))\wh (y))-\al\log (\J^\be(\nu))\wh (y)
\end{equation}
\end{wn}

\begin{stw}
For $\be>0$  and $\nu\in ID$ we have
\begin{equation}
\mathcal{J}^{\be}(\mathcal{J}^{\be}(\nu))=
\mathcal{L}(\int_0^1t\,dY_{\nu}( r_{(\be,\be)}(t))),
\end{equation}
where $r_{(\be,\be)}(u)\,:=\,u^{\be}(1-\beta\,\log u)$ is an
increasing time change in a L\'evy process $Y$.
\end{stw}

\begin{rem}
Let us note that for the functions $r_{(\al,\be)}, \ \ \al\neq\be$,
and $r_{(\be,\be)}$, from Propositions 1 and 2, we have
\[
\lim_{\al\to \be}r_{(\al,\be)}(u)= r_{(\be,\be)}(u), \ \ \
\mbox{for}\ 0\le u \le 1 .
\]
\end{rem}
\begin{rem}
For $\be=1$, probability measures from
$\mathcal{U}^{<2>}:=\J^1(\J^1(ID))$ were called \emph{2-times
s-selfdecomposable  distributions} in Jurek (2004). In fact,
\emph{m-times s-selfdecomposability} was defined there
inductively, and the corresponding classes $\mathcal{U}^{<m>}$,
for $m=1,2,...$\, were described in many ways; cf. Propositions 3
and 4.
\end{rem}

Each infinitely divisible probability distribution $\mu$ is
uniquely determined by  a triplet: a shift vector $a$, a Gaussian
covariance operator $R$ and a L\'evy spectral measure $M$ that
appear in the L\'evy-Khintchine formula, as it is recalled at
beginning of Section 4. Therefore,  following the notation from
Parthasarathy (1967), Chapter VI, we will write that
$\mu=[a,R,M]$.

Directly from Jurek (1988), or from Lemma 1 below, if $\mu=[a,R,M]$
and $\J^\be(\mu)=[a^{(\be)}, R^{(\be)}, M^{(\be)}]$  and
\begin{equation}
b_{M,\be}:=\int_{\{||x||>1\}} x\,||x||^{-1-\be}\,M(dx)\in \Rset^d \
( \ \mbox{or} \  E )
\end{equation}
then we have
\begin{multline}
a^{(\be)}:= \be(\be+1)^{-1}(a+ b_{M, \be}), \ \ \   R^{(\be)}:=\be(2+\be)^{-1}R  \\
M^{(\be)}(A):=\int_0^1\,T_{t^{1/\be}}M(A)dt=\int_0^1\int_{\Rset^d}1_{A}(t^{1/\be}\,x)M(dx)dt,
\ \  \mbox{for} \ A\in \mathcal{B}_0.
\end{multline}
The above $\mathcal{B}_0$ stands for all Borel subsets of
$\Rset^d\setminus\{0\}$ (or $E\setminus\{0\}$ if one consider
results on Banach space $E$). Note that  one needs to change the
order of integration in the formula (1.10) in Jurek (1988) to get
the above form of $a^{(\be)}$.

\begin{stw}
For positive $\al\neq\be$, if $\mu=[a,R,M]$ and
$\J^\al(\J^\be(\mu))=[a^{(\al,\be)},R^{(\al,\be)},M^{(\al,\be)}]$
then
\begin{description}
\item[(i)] $a^{(\al,\be)}=\al \be(1+\al)^{-1}(1+\be)^{-1}\ a +\al
\be(\be-\al)^{-1}\\ \Big[(\al+1)^{-1}\int_{\{||x||>1\}} x
||x||^{-1-\al}M(dx) -
(\be+1)^{-1}\int_{\{||x||>1\}}x ||x||^{-1-\be} M(dx)\Big]\\
=
\frac{\be}{\be-\al}\,a^{(\al)}-\frac{\al}{\be-\al}\,a^{(\be)}$\,.

\item[(ii)]$
R^{(\al,\be)}=\tfrac{\al}{2+\al}\cdot\tfrac{\be}{2+\be}R =
\frac{\be}{\be-\al}R^{(\al)}- \frac{\al}{\be-\al}R^{(\be)}$\,.

\item[(iii)]${\displaystyle
M^{(\al,\be)}(A)=\int^{1}_{0}\int^{1}_{0}T_{t^{1/\al}s^{1/\be}}M(A)\;ds\;dt}
=\frac{\be}{\be-\al}\,M^{(\al)}(A)-\frac{\al}{\be-\al}\,M^{(\be)}(A)$\,
\end{description}
for  all Borel sets $A$ in $\mathcal{B}_0$.
\end{stw}

\begin{rem} Let us remark here that the  parameters in the triple
corresponding to the composition $\J^\al\circ\J^\be$ are linear
combinations of parameters corresponding to the mappings $\J^\al$
and $\J^\be$ in an identical way as they appear in the formula of
$r_{(\al,\be)}$ in Proposition 1.
\end{rem}

\medskip
From Jurek(1988), Corollary 1.1, if  $0<\al<\be$ then we have
$\mathcal{U}_{\al}\subseteq \mathcal{U}_{\be}$. The converse
inclusion is determined as follows.
\begin{wn}
Let $0<\al<\be$. In order that
$\J^{\be}\bigl(\rho\bigr)\in\mathcal{U}_{\al}$, for some $\rho\in
ID$, it is necessary and sufficient that  $\rho$ admits the
following convolution factorization $\rho =\J^\alpha\bigl(\nu^{
* (1-\alpha/\beta)}\bigr)*\nu^{
* \alpha/\beta}=\big(\J^\al(\nu)\big)^{\ast(1- \al/\be)}\ast \nu^{\ast \al/\be}$,
for some $\nu\in \textit{ID}.$
\end{wn}

\medskip
We will illustrate our results in the following examples.

\textbf{3. Examples.}

\medskip
\textbf{(a)} For $\al>0$ and $\be:=2\al$, Proposition 1 gives that
$r_{(\al,2\al)}(u)=2u^{\al}-u^{2\al}$ and
$r_{(\al,2\al)}^{-1}(t)=\big(1-\sqrt{1-t}\big)^{1/\al}$, and we have
that:
\begin{multline} \J^\al\bigl(\J^{2\al}\bigl(\nu\bigr)\bigr)= \La \big(\int^{1}_{0}\,u
\,dY_{\nu}(2u^\al- u^{2\al})\big)= \La\ \big(\int^{1}_{0}\Bl
1-\sqrt{1-t}\Br^{1/\al}dY_{\nu}(t)\big)= \\
\La\ \big(\int^{1}_{0}\Bl
1-\sqrt{1-t}\Br^{1/\al}d\widetilde{Y}_{\nu}(t)\big)=
\La\Bigl(\int^{1}_{0}\Bl 1-\sqrt{t}\Br^{1/\al}dY_{\nu}(t)\Bigr).
\end{multline}
where $\widetilde{Y}_{\nu}(t):=Y_{\nu}(1)-Y_{\nu}(1-t), 0\le t\le
1$, and $Y_{\nu}(t), 0\le t\le 1$, have the same distributions, and
$d\widetilde{Y}_{\nu}(t)=\,d\,Y_{\nu}(t)$.

\medskip
\begin{rem}
From Proposition 3 we may get the formulae for the triple
$a^{(\al,2\al)}$, $R^{(\al,2\al)}$ and $M^{(\al,2\al)}$. However,
since in our example we have an explicit form for
$r_{(\al,2\al)}^{-1}$ therefore applying Lemma 1 and the formula
(13) to the last integral in (9) we have an alternative way of
getting the triplet in question. Thus we have
\[
M^{(\al,2\al)}(A)=\int_0^1T_{(1-\sqrt{t})^{1/\al}}M(A)dt=
\int_0^1\int_{\Rset^d}\,1_A((1-\sqrt{t})^{1/\al}\,x)M(dx)\,dt\,,
\]
 for all Borel sets $A$ in $\mathcal{B}_0$.
\end{rem}

\medskip
\textbf{(b)} Let $\sigma_p:=[a,0, M_p]$ denotes \emph{a stable
distribution} with an exponent $0<p <2$, that is, for some finite
measure $\gamma$ on the unit sphere $S=\{x: ||x||=1\}$ we have
\begin{equation}
M_p(A):=\int_S\int_0^{\infty}1_A(r\,x)r^{-p-1}dr\gamma(dx), \ \ A\in
\mathcal{B}_0 \,;
\end{equation}
cf. Araujo-Gine (1980), Chapter 3, Theorem 6.15.  Then
\[
b_{M_p,\be} = (\be+p)^{-1} \bar{\gamma}, \ \ \mbox{where} \ \
\bar{\gamma}:=\int_Su\,\gamma(du) ; \ \
M_p^{(\be)}=\be(\be+p)^{-1}\,M_p,
\]
by (7) and (8), respectively. Consequently, from Proposition 3 we
get
\begin{multline}
\J^\al(\J^\be(\sigma_p))\\
=\Big[\,\frac{\al\be}{(\al+1)(\be+1)}\, \big( a
+\frac{\al+\be+p+1}{(\al+p)(\be+p)}\,\,\bar{\gamma}\big), \ \ 0, \ \
\frac{\al\be}{(\al+p)(\be+p)}\,M_p\,\Big]\\
= \sigma_p^{\ast\, \frac{\al\be}{(\al+p)(\be+p)}} \ast \delta_{x_0},
\qquad \qquad \qquad \qquad \qquad
\end{multline}
where a vector $x_0$ is given by the formula
\[
x_0:=\frac{\al
\be\,(\al+\be+p+1)}{(\al+1)(\be+1)(\al+p)(\be+p)}\,\big
[(p-1)\,a+\bar{\gamma}\,\big],
\]
that is, in (11), up to a shift vector, we get  a convolution power
of the stable measure $\sigma_p$.

\medskip
\textbf{(c)} Let $e_{\lambda}$ denotes the exponential distribution
with the parameter $\lambda$. Then its L\'evy spectral measure $M_e$
has the density $e^{-\lambda\,x}\,x^{-1}1_{(0, \infty)}(x)$  and has
the characteristic function
\[
\frac{1}{\lambda-i\,y}=\widehat{e_{\lambda}}(y)=\exp\int_0^{\infty}(e^{iyx}-1)\frac{e^{-\lambda\,x}}{x}
dx\, , y\in\Rset\,.
\]
Then the L\'evy spectral measure $M_e^{(\be)}$  has the density
$\lambda\be(\lambda\,x)^{\be-1}\Gamma(-\be,\lambda\,x)$, $x>0$,
where for $c\in\Rset$
\[
\Gamma(c,x)=:\int_x^{\infty}u^{c-1}e^{-\lambda\,u}du , x>0, \ \
\mbox{is the incomplete Euler gamma function};
\]
cf. Gradsteyn and Ryzhik (1994), Section 8.3  for other
representations of the gamma function. Then  for $A\subset
(0,\infty)$, we get the formula
\begin{equation}
M_e^{(\al,\be)}(A)= \frac{\al\be\lambda}{\be-\al}\int_A\,
[(\lambda\,x)^{\al-1}\Gamma(-\al,\lambda\,x)-(\lambda\,x)^{\be-1}\Gamma(-\be,\lambda\,x)]\,dx,
\end{equation}
i.e., the density is a combinations (with variable coefficients) of
two incomplete gamma (Euler) functions.

\medskip
\textbf{4. Proofs and auxiliary results.} Let us recall that for a
probability Borel measures $\mu$ on $\Rset^d$, its
\emph{characteristic function} $\hat{\mu}$ is defined as
\[
\widehat{\mu}(y):=\int_{\Rset^d} e^{i<y,x>}\mu(dx), \ y\in\Rset^d,
\]
where $<\cdot,\cdot>$ denotes the scalar product (or a bilinear form
between the conjugate Banach space $E^{\prime}$ and the space $E$).
Recall that for infinitely divisible measures $\mu$ their
characteristic functions admit the following L\'evy-Khintchine
formula
\begin{multline}
\widehat{\mu}(y)= e^{\Phi(y)}, \ y \in \Rset^d, \ \ \mbox{and the
L\'evy exponents} \ \ \Phi(y)=i<y,a>- \\ \frac{1}{2}<y,Ry> +
\int_{\Rset^d \backslash \{0 \}}[e^{i<y,x>}-1-i<y,x>1_B(x)]M(dx),
\qquad \qquad \qquad
\end{multline}
where $a$ is  a \emph{shift vector}, $R$ is a \emph{covariance
operator} corresponding to the Gaussian part of $\mu$ and $M$ is a
\emph{L\'evy spectral measure}. Since there is a one-to-one
correspondence between a measure $\mu \in ID$ and the triples $a$,
$R$ and $M$ in its L\'evy-Khintchine formula (13) we will write
$\mu=[a,R,M]$. Finally, let recall that\footnote{The integrability
criterium (14) is true also in real separable Hilbert spaces, cf.
Parthasarathy (1967), Chapter VI. But no such a characterization is
available for infinite dimensional Banach spaces; cf. Araujo-Gine
(1980).} \emph{
\begin{equation}
M \ \mbox{is L\'evy spectral measure on $\Rset^d$ iff} \ \
\int_{\Rset^d}\min(1, ||x||^2)M(dx)<\infty
\end{equation}
} For this note it is important to recall the following crucial
fact.
\begin{lem}
If the random integral $X:= \int_I h(t)dY(r(t))$ exists then we
have
\begin{equation}
\log (\mathcal{L}(X))^{\wh} (y)=\int_I \log
(\mathcal{L}(Y(1))^{\wh}(h(s)y)dr(s)=\int_I\,\Phi(h(s)y)dr(s), \ \
y\in \Rset^d
\end{equation}
where $\Phi$ is the L\'evy exponent of $(\mathcal{L}(Y(1)))^{\wh}$.
\end{lem}

The formula is a straightforward consequence of our definition
(integration by parts)  of the random integrals (1). The proof is
analogous to that in  Jurek-Vervaat(1983), Lemma 1.1 or
Jurek-Mason (1993), Lemma 3.6.4 or Jurek (1988), Lemma 2.2 (b).

Note that for  bounded sets $I\subset\Rset^+$ and continuous $h$
integrals of the form (1)  are well-defined. In particular, we
have
\begin{equation}
(\J^\beta(\nu))^{\wh}(y) =\exp\int^{1}_{0} \log
\hat{\nu}(t^{1/\be}y)\,dt, \ \ \ y \in \Rset^d \ \ (\mbox{or} \ \ y
\in E^{\prime}),
\end{equation}

Since $\mathcal{J}^{\beta}$ are mappings from $ID$ into $ID$,
therefore we may consider their compositions. Here we recall their
basic properties for further references.
\begin{lem} \emph{Relations between the
mappings $\mathcal{J^{\beta}}$, for $\be>0$.}

(a) Each mapping $\mathcal{J^{\beta}}$ is a continuous  isomorphism
between convolution semigroups $ID$ and $\mathcal{U}_{\beta}$.

(b) For $\be, c>0$ and $\nu\in ID$ we have
\begin{equation}
(\mathcal{J^{\beta}}(\nu))^{\ast c}=\mathcal{J^{\beta}}(\nu^{\ast
c})
\end{equation}

(c)\,The composition of the mappings $\mathcal{J^{\be}}, \be >0,$ is
commutative, that is, for
 $\al, \be > 0$ and $\nu\in ID$
\begin{equation}
\J^\al\big(\J^\be(\nu)\big)= \J^\be\big(\J^\al(\nu)\big). \ \ \
\end{equation}
\end{lem}
\emph{Proof of Lemma 2.} Parts (a) and (b) have  proofs along the
lines of a proof of Theorem 1.3 (a) , (b) and (c) in Jurek (1988)
for the operator $Q=I$. But basically we utilize the formula (16).
For the commutativity property note that using the formula (16) we
get
\begin{align}
\lefteqn{\Bl\J^\alpha\left(\J^\beta\left(\nu\right)\right)\Br\wh\left(y\right)=}\nn\\
&\qquad=\exp{\Bl\int^{1}_{0}\log{ \bl\J^{\beta}\left(\nu\right)\br\wh \left(t^{1/\alpha}y\right)}\,dt\Br}
=\exp{\Bl\int^{1}_{0} \int^{1}_{0}\log{\hat{\nu} \left(s^{1/\beta}t^{1/\alpha}y\right)}\,ds\,dt\Br}\nn\\
&\qquad=\exp{\Bl\int^{1}_{0}\int^{1}_{0}\log{\hat{\nu} \left(s^{1/\beta}t^{1/\alpha}y\right)}\,dt\,ds\Br}
=\J^\beta\bl\J^\alpha\left(\nu\right)\br\wh\left(y\right),\nn
\end{align}
which indeed proves the commutativity in (c). This completes a proof
of Lemma 2.

\medskip
\emph{Proof of Proposition 1.} Note that the above can be rewritten
as follows
\begin{multline}
\log{\bl\J^\alpha(\J^\beta(\nu))\wh (y)}
=\int^{1}_{0}\int^{1}_{0}\log{\hat{\nu}
\left(s^{1/\beta}t^{1/\alpha}y\right)}\,dt\,ds \ \ \
(\mbox{putting}\ \ t:= u^{\al}s^{-\al/\be}) \\
=\int^{1}_{0}\int^{s^{1/\beta}}_{0}\log{\hat{\nu} \left(
uy\right)}\frac{\alpha u^{\alpha-1}}{s^{\alpha/\beta}}\,du\,ds
=\int^{1}_{0}\log{\hat{\nu}(uy)}\alpha u^{\alpha-1}
\int^{1}_{u^\beta}s^{-\alpha/\beta}\,ds\,du \\
= \frac{\beta}{\beta-\alpha}\int^{1}_{0}\log{\hat{\nu }(uy)}\alpha
u^{\alpha-1}\,du -
\frac{\alpha}{\beta-\alpha}\int^{1}_{0}\log{\hat{\nu}(uy)}\beta
u^{\beta-1}\,du \\
=\int_0^1\log\hat{\nu}(uy) d\, r_{(\al,\be)}(u)= \log\Big(
\La(\int^{1}_{0}u\,dY_{\nu}(r_{(\al,\be)}(u)\br\bigr)^{\wh}(y)
\big),
\end{multline}
where the function $r_{(\al,\be)}$ is given in Proposition 1 and the
last equality follows from Lemma 1 .

\medskip
\emph{Proof of Corollary 1.} From the property (17) and then (19) we
get
\begin{multline}
\log{\bl\J^\alpha(\J^\beta(\nu^{\ast (\be-\al)}))\wh (y)}=(\be-\al) \log{\bl\J^\alpha(\J^\beta(\nu))\wh (y)} \\
= \beta\int^{1}_{0}\log{\hat{\nu}(uy)}\alpha u^{\alpha-1}\,du -
\alpha\int^{1}_{0}\log{\hat{\nu}(uy)}\beta u^{\beta-1}\,du \\
=\beta\int^{1}_{0}\log{\hat{\nu }( t^{1/\al}y)}\,dt -
 \alpha\int^{1}_{0}\log{\hat{\nu }\left(
 t^{1/\be}y\right)}\,dt \ \ \ \ \ \\
 =\int^{1}_{0}\log(\nu^{\ast \beta})\wh(t^{1/\al}y)\,dt -
 \int^{1}_{0}\log (\nu^{\ast\alpha})\wh
 (t^{1/\be}y)dt \ \ \ \  (\mbox{by (17)}) \\
= \log(\J^\alpha(\nu^{\beta})\wh(y)- \log(\J^\beta( \nu^{\ast
\alpha}))\wh(y)).
\end{multline}

In other words we have
\begin{equation}
(\J^\alpha(\J^\beta(\nu^{\ast(\be-\al)}\bigr)))\wh\left(y\right)\cdot
(\J^\beta\bigl(\nu^{\ast \alpha} \bigr))\wh(y)
=(\J^\alpha(\nu^{\ast\beta}))\wh(y)
\end{equation}
But (21), in terms of probability measures, coincides with the
formula (4) that completes the proof of the Corollary.

\medskip
\emph{Proof of Proposition 2.} Similarly as at the beginning of
(15) we get
\begin{multline}
\log{\bl\J^\beta(\J^\beta(\nu))\wh (y)}
=\int^{1}_{0}\int^{1}_{0}\log{\hat{\nu}
\left(s^{1/\beta}t^{1/\beta}y\right)}\,dt\,ds \\
=\int^{1}_{0}\int^{s^{1/\beta}}_{0}\log{\hat{\nu}} \left(
uy\right)\frac{\be u^{\be-1}}{s}\,du\,ds =
\int^{1}_{0}\log{\hat{\nu}}\left(uy\right) \be u^{\be-1} (-\be
\log u)du \\
=
\int^{1}_{0}\log{\hat{\nu}}\left(uy\right)d[u^{\be}(1-\be\log
u)]du,
\end{multline}

which, via Lemma 1, is the  the statement (6) in Proposition 2.

\medskip
\emph{Proof of Proposition 3.} First, note that for $M^{(\be)}$
given by (8), using (7) we get
\begin{multline*} b_{M^{(\be)},\al}
=\int_0^1\int_{\Rset^d}1_{B^c}(x)\,x\,||x||^{-1-\al}M(t^{-1/\be}dx)dt\\
=\int_0^1 \int_{\Rset^d}
1_{B^c}(t^{1/\be}x)\,x\,||x||^{-1-\al}t^{-\al/\be}M(dx)dt\\=
\int_{\{||x||>1\}}\,x\,||x||^{-1-\al}\int_{||x||^{-\be}}^1\,t^{-\al/\be}dt\,M(dx)=\be(\be-\al)^{-1}(b_{M,\al}-b_{M,\be}).
\end{multline*}

Second,  successively using (8)  (for the shift vector) and the
above one gets
\begin{multline*}
(a^{(\be)})^{(\al)}= \al(\al+1)^{-1}[a^{(\be)}
+ b_{M^{(\be)},\al}] \\
=\al(\al+1)^{-1}[\be(\be+1)^{-1}a + \be(\be+1)^{-1} b_{M,\be} +
\be(\be-\al)^{-1}(b_{M,\al}-b_{M,\be})]\\
=\al \be(1+\al)^{-1}(1+\be)^{-1}\ a +\al \be(\be-\al)^{-1}
[(\al+1)^{-1}b_{M,\al} -(\be+1)^{-1}b_{M,\be}]\\
= \frac{\be}{\be-\al}\,a^{(\al)}-\frac{\al}{\be-\al}\,a^{(\be)}, \ \
\ \ \ \ \ \
\end{multline*}
and the last equality one checks by straightforward computation.
This proves the formula for the shift vector \textbf{(i)}. Part
\textbf{(ii)} follows easily from (5). Finally we have
\begin{multline*}
M^{(\al,\be)}(A)=
(M^{(\al)})^{(\be)}(A)=\int_0^1T_{t^{1/\be}}M^{(\al)}(A)dt\\
=\int_0^1\int_0^1\,T_{t^{1/\be}s^{1/\al}}M(A)dtds=\int_0^1\int_0^{s^{\be/\al}}\,T_{w^{1/\be}}M(A)s^{-\be/\al}dw\,ds\\
=\int_0^1\int_{w^{\al/\be}}^1\,T_{w^{1/\be}}M(A)s^{-\be/\al}ds\,dw=
\al(\al-\be)^{-1}\int_0^1(1-w^{(\al-\be)/\be})T_{w^{1/\be}}(A)dw \\
=\al(\al-\be)^{-1}M^{(\be)}(A)-\al(\al-\be)^{-1}\int_0^1T_{w^{1/\be}}(A)w^{\al/\be-1}dw
\\ =\al(\al-\be)^{-1}M^{(\be)}(A)-\be(\al-\be)^{-1}M^{(\al)}(A) = \be(\be-\al)^{-1}M^{(\al)}-\al(\be-\al)^{-1}M^{(\be)},
\end{multline*}
which gives  the equality \textbf{(iii)}. Thus the proof of
Proposition 3 is completed.

\emph{Proof of Corollary 2.}  Replacing $\nu$ by $\nu^{\ast 1/\be}$
in Corollary 1 and using the commutativity of the mappings
$\J^{\al}$ and $\J^{\be}$ (Lemma 2 c)), we get
\begin{equation}
\J^\be\bigl(\J^\alpha\bigl(\nu^{\ast(1 -\alpha/\beta)}\bigr)*
\nu^{\ast \alpha/\beta}\bigr)= \J^\al\bigl(\nu \bigr),
\end{equation}
and hence for $\rho :=\J^\alpha\bigl(\nu^{(1-\alpha/\beta)}\bigr)*
\nu^{\ast\, \alpha/\beta}$ we get that
$\J^{\be}(\rho)\in\mathcal{U}_{\al}$. (Note that from the above we
 may also conclude that $\mathcal{U}_{\al}\subseteq
\mathcal{U}_{\be}$).

\noindent Conversely, let $\J^{\be}\bigl(\rho\bigr)=\J^\al\bigl(\nu
\bigr)$, for some $\nu\in ID$. Hence the equality (23) implies that
$\rho =\J^\alpha \bigl(\nu^{\ast(1- \alpha/\beta)}\bigr)\ast
\nu^{\ast \alpha/\beta}$, because $\J^{\be}$ is one-to-one mapping
(Lemma 2 a)). Thus this completes the proof.

\medskip
\medskip
\textbf{Acknowledgement.} This research  has been in a part
supported also by a Maria Curie Transfer of Knowledge (TOK)
Fellowship of the European Community's Sixth Framework under
contract number MTKD-CT-2004-013389. Both Authors would like to
thank Professor Ewa Damek, chairwoman of the TOK program, for her
help.

\begin{center}
\textbf{References}
\end{center}

\medskip
\noindent [1] A. Araujo and E. Gine (1980). \emph{The central limit
theorem for real and Banach valued random variables.} John Wiley \&
Sons, New York.

 \noindent [2] T. Aoyama and M. Maejima (2007).
Characterizations of subclasses otf type G distributions on
$\Rset^d$ by stochastic random integral representation,
\emph{Bernoulli}, vol. 13, pp. 148-160.

\noindent[3] A. Czyżewska-Jankowska, Z. J. Jurek (2008). Random
integral representation of the class $L^f$ distributions and some
related properties; \emph{submitted}.

\noindent[4] I. S. Gradshteyn and I. M. Ryzhik (1994). \emph{Table
of integrals, series, and products}, 5th Edition, Academic Press.

\noindent [5] Z. J. Jurek (1996). s-stability and random integral
representations of limit laws, \emph{Demonstration Math.}, vol. XXIX
, No 2 , pp. 363-370.

\noindent [6] Z. J. Jurek (2004). The random integral representation
hypothesis revisited: new classes of s-selfdecomposable laws. In:
Abstract and Applied Analysis; \emph{Proc. International Conf.
ICAAA} , Hanoi, August 2002, str. 495-514. World Scientific,
Hongkong.

\noindent [7] Z. J. Jurek (1988). Random Integral representation for
Classes of Limit Distributions Similar to Levy Class $L_{0}$,
 \emph{Probab. Th. Fields.} 78, pp. 473-490.

\noindent [8] Z. J. Jurek and W. Vervaat (1983). An integral
representation for selfdecomposable Banach space valued random
variables, \emph{Z. Wahrscheinlichkeitstheorie verw. Gebiete}, 62,
pp. 247-262.

\noindent[9] K. R. Parthasarathy (1967). \emph{Probability measures
on metric spaces}. Academic Press, New York and London.

\medskip
\noindent
Institute of Mathematics \\
University of Wroc\l aw \\
Pl.Grunwaldzki 2/4 \\
50-384 Wroc\l aw, Poland \\
e-mail: zjjurek@math.uni.wroc.pl \ \ or \ \ czyzew@math.uni.wroc.pl \\
www.math.uni.wroc.pl/$^{\sim}$zjjurek

\end{document}